\documentclass[final]{siamltex}

\usepackage{amssymb,amsmath}
\newcommand{\beq}       {\begin{equation}}
\newcommand{\beqstar}       {\begin{equation*}}
\newcommand{\eeq}       {\end{equation}}
\newcommand{\eeqstar}       {\end{equation*}}
\newcommand{\bery}      {\begin{array}}
\newcommand{\eery}      {\end{array}}
\newcommand{\berys}     {\begin{array*}}
\newcommand{\eerys}     {\end{array*}}
\newcommand{\beqry}     {\begin{eqnarray}}
\newcommand{\eeqry}     {\end{eqnarray}}
\newcommand{\beqrys}    {\begin{eqnarray*}}
\newcommand{\eeqrys}    {\end{eqnarray*}}

\newcommand{\R}{{\mathbf R}}        

\def\myref#1{(\ref{#1})}

\def\Sh{\mathop{\rm Sh}\nolimits}

\def\Sh{\mathop{\rm Sh}\nolimits}

\def\diag{\mathop{\rm diag}\nolimits}
\def\GL{\mathop{\rm GL}\nolimits}
\def\Gr{\mathop{\rm Gr}\nolimits}
\def\R{{\mathbf R}}

\begin{document}

\title{Birth of the shape of a reachable set\thanks{This
        work was partially supported by the Russian Foundation for Basic Research (grants  11-08-00435, 13-08-00441).}}


\author{Elena Goncharova\thanks{Institute for System Dynamics and Control Theory, Siberian Branch of Russian Academy of Sciences,
        134, Lermontov st., 664033 Irkutsk, Russia
        ({\tt goncha@icc.ru}).}
        \and Alexander Ovseevich\thanks{Institute for Problems in Mechanics, Russian Academy of Sciences, 101/1, Vernadsky av.,
        119526 Moscow, Russia ({\tt ovseev@ipmnet.ru}).}}

\maketitle
\begin{abstract}
We study reachable sets of linear dynamical systems with geometric constraints
on control. The main result is the existence of a limit shape of the reachable
sets as the terminal time tends to zero. Here, a shape of a set stands for the
set regarded up to an invertible linear transformation.
\end{abstract}

\begin{keywords}
Linear dynamical system, reachable sets.
\end{keywords}

\begin{AMS}
93B03, 34H05
\end{AMS}

\pagestyle{myheadings} \thispagestyle{plain} \markboth{E. GONCHAROVA AND A. OVSEEVICH}{BIRTH OF THE
SHAPE}

\noindent
\section{Problem statement}
Consider the following linear 
control system:
\begin{equation}\label{syst}
\dot x=Ax+Bu,\; x\in \mathbb{V}=\R^n,\; u\in U\subset \mathbb{U}=\R^m,
\end{equation}
where $U$ is a central symmetric convex  body in $\R^m$. Given $T>0$,
we get the reachable set ${\cal D}(T)$ being the set of the ends $\{x(T)\}$ at
time $T$ of all admissible trajectories of system \myref{syst} such that
$x(0)=0.$ We study the asymptotic behavior of the reachable sets ${\cal D}(T)$
as $T\to 0$. The problem is not as trivial as it may seem, since the obvious
answer ${\cal D}(T)\to\{0\}$ is not good enough. Of course, ${\cal D}(T)$ is
small with a small $T$. Still we can look at it through a kind of microscope
and see fine details of its development, in particular, its shape, where a
shape of a set is the set regarded up to an invertible linear transformation.
For instance, consider a rather trivial system \myref{syst}, where $A=0$,
$B=1$, and $U$ is the unit ball $B_1$ in $\R^n=\R^m$. Then ${\cal D}(T)=TB_1$,
and the reachable set has the shape of a ball at any time $T>0$.
In our example,  for all $T>0$,
the sets ${\cal D}(T)$ have the same shape, so that
there exists a limit of the shapes as $T\rightarrow  0$, and this limit is not zero, but the shape of a ball. 

In this paper we will show that, for any completely controllable linear system, shapes of the
reachable sets converge as time tends to zero. We will give an estimate for the rate of convergence.
The paper can be regarded as the $T \to 0$ counterpart of the results \cite{ovs} on the large-time
dynamics of the reachable sets, where the existence of a limit shape
has been established.

\section{Shapes of reachable sets}\label{shapes section}

There is a well developed mathematical concept of shapes, see \cite{ovs}.
Consider the metric space $\mathbb{B}$ of central symmetric convex bodies with
the Banach-Mazur distance $\rho $:
\begin{eqnarray*}\label{B-M}
\rho (\Omega _1,\Omega _2)=\log (t(\Omega _1,\Omega _2)t(\Omega _2,\Omega _1)), \; t(\Omega_1 ,
\Omega_2 ) = \inf \{t \geq 1: t\Omega_1\supset \Omega_2 \}.
\end{eqnarray*}
The general linear group $\GL(\mathbb{V})$ naturally acts on the space
$\mathbb{B}$ by isometries. The factorspace $\mathbb{S}$ is called the space of
shapes of central symmetric convex bodies, where the shape $\Sh\Omega \in
\mathbb{S}$ of a convex body $\Omega \in \mathbb{B}$ is the orbit $\Sh\Omega =
\{g \Omega : {\rm det}\, g \neq 0\}$ of the point $\Omega$ with respect to the
action of $\GL(\mathbb{V})$. The Banach-Mazur factormetric $$\rho (\Sh\Omega
_1,\Sh\Omega _2)=\inf_{g\in GL(\mathbb{V})}\rho (g\Omega _1,\Omega _2)$$ makes
$\mathbb{S}$ into a compact metric space.   If $\rho (\Omega _1(T),\Omega _2(T)) \to 0$ as $T\to 0$, 
the convex bodies $\Omega_1$ and $\Omega_2$ are asymptotically equal, and we
write $\Omega _1(T) \sim \Omega _2(T)$. The asymptotic equivalence of shapes is
defined similarly.

\section{Main result: Autonomous case}


We assume that system \myref{syst} is time-invariant and the Kalman
controllability condition holds. The Kalman condition ensures that the
reachable sets ${\cal D} (T)$ to system \myref{syst} are central symmetric
convex bodies in $\R^n$. In what follows, the convergence of 
shapes of reachable sets is understood in the sense of the Banach-Mazur metric.
\begin{theorem}\label{thmain} The shapes $\Sh{\cal D} (T)$
have a limit $\Sh_0 $ as $T\to 0$. The Banach--Mazur distance $\rho(\Sh{\cal D} (T),\Sh_0 )$ is
$O(T)$.
\end{theorem}

\noindent This means that there exists a time independent convex body $\Omega$ such that 
$${\cal D} (T)\sim C(T)\Omega,$$ where $C(T)$ is a matrix function. The Banach-Mazur distance between the left- and
the right-hand sides of the latter formula is $O(T)$.


The proof of Theorem~\ref{thmain} is based on two easy lemmas, and the use of
the Brunovsky normal form of a controllable system.
\begin{lemma} \label{linear_feedback}
Consider the linear system
\begin{equation}\label{syst2}
\dot x={\widetilde A}x+{\widetilde B}u, 
\; u\in U,
\end{equation}
obtained from system {\rm \myref{syst}} by adding a linear feedback, that is, ${\widetilde A}=A+BC$
and ${\widetilde B}=B$. Then the Banach--Mazur distance $\rho({\cal D} (T),\widetilde{\cal D} (T))$
is $O({T})$ as $T\to0$, where ${\cal D} (T)$ and $\widetilde{\cal D} (T)$ are the reachable sets to
systems {\rm \myref{syst}} and {\rm \myref{syst2}}, respectively.
\end{lemma}

\begin{lemma} \label{gauge}
Consider the linear system
\begin{equation}\label{syst22}
\dot x={\widetilde A}x+{\widetilde B}u, 
\; u\in U, 
\end{equation}
obtained from system {\rm \myref{syst}} by a gauge transformation, where ${\widetilde A}=C^{-1}AC$,
${\widetilde B}=C^{-1}B$, and $C$ is an invertible matrix.  Then $\Sh\widetilde{\cal D} (T)=\Sh{\cal
D} (T)$, where $\widetilde{\cal D} (T)$ is the reachable set to system {\rm \myref{syst22}}.
\end{lemma}

From the Lemmas it follows that applying gauge transformations coupled with
adding a linear feedback do not affect the validity of Theorem~\ref{thmain}.
However, by means of these transformations, a general system \myref{syst} can
be reduced to the Brunovsky normal form~\cite{Brun}. One can show that shapes $\Sh{\cal D} (T)$ of reachable sets of the
Brunovsky system do not depend on $T$.

The Brunovsky classification is known to be closely related to the Grothendieck
theorem \cite{Groth} on decomposition of vector bundles on ${\mathbf P}^1$ into
a sum  of linear bundles.

\section{Non-autonomous case} 

Consider system \myref{syst}, where now the data $A$, $B$, and $U$ are
$C^\infty$-functions of  time $t\geq0$. By a standard trick we make
\myref{syst} into the time-invariant system
\begin{eqnarray*}\label{timeinv}
&&\dot\tau=1,\\ &&\dot x={ A(\tau)}x+{ B(\tau)}u.
\end{eqnarray*}
Consider the Lie algebra ${\mathcal L}$ generated by the vector fields $\left(
1,
  { A(\tau)}x \right)$ and $\left( 0,
  { B(\tau)}u \right)$ in $\R\times\mathbb{V}=\R^{n+1}$, where $u\in\R^m$ is a constant vector. Define
  ${\mathcal L}(\tau,x)$ as the set of the values at $(\tau,x)$ of all vector fields from ${\mathcal L}$.

Assume that the following Kalman type condition holds:
\begin{quote}
For each $(\tau,x)\in\R\times\mathbb{V}$ the set ${\mathcal L}(\tau,x)$ coincides with the entire
tangent space $\R^{n+1}$. In other words,
\begin{equation}\label{ext_kalman}
\dim{\mathcal L}(\tau,x)=n+1.
\end{equation}
\end{quote}

\begin{theorem}\label{thmain2} Let
$\mathcal{D}(T)$ be the reachable set to a non-autonomous system of the form {\rm \myref{syst}} and
the genericity condition {\rm \myref{ext_kalman}} hold. The shapes $\Sh{\cal D} (T)$ have a limit
$\Sh_0 $ as $T\to0$. Moreover, the Banach--Mazur distance $\rho(\Sh{\cal D} (T),\Sh_0 )$ is
$O({T}).$
\end{theorem}

The proof is based on a reduction to the case $A=0$ and a subsequent study of
the filtration $F_k^*=\{\xi \in \mathbb{V}^*:{ B(t)}^*\xi=O(t^k)\}$ of the dual
space $\mathbb{V}^*$.


\begin{thebibliography}{9}

\bibitem{ovs}  {\it Ovseevich~A.I.} Asymptotic behavior of attainable and
superattainable sets // Proceedings of the Conference on Modeling, Estimation
and Filtering of Systems with Uncertainty, Sopron, Hungary, 1990,
Birkha\"{u}ser, Basel, Switzerland. 1991. P. 324--333.

\bibitem{fig-periodic} {\it Figurina T.Yu., Ovseevich A.I.}  Asymptotic behavior
of attainable sets of linear periodic control systems
// J. Optim. Theory Appl. 1999. Vol. 100. No 2. P. 349--364.

\bibitem{Brun} {\it Brunovsky~P.} A classification of linear controllable systems. Kibernetika 6 (1970), 176-188.

\bibitem{Groth} {\it Grothendieck A.} Sur la classification des fibr\'es holomorphes sur la sph\'ere de Riemann,  Amer. J. Math., 79, 121-138 (1957)




\end{thebibliography}
\end{document}